\documentclass{amsart}

 \newtheorem{thm}{Theorem}[section]
 \newtheorem{cor}[thm]{Corollary}
 \newtheorem{lem}[thm]{Lemma}
 \newtheorem{prop}[thm]{Proposition}
 \theoremstyle{definition}
 \newtheorem{defn}[thm]{Definition}
 \theoremstyle{remark}
 
 \numberwithin{equation}{section}

 \newcommand{\norm}[1]{\left\Vert#1\right\Vert}

\begin{document}

\title[Modular Lattice for $C_{o}$-Operators.]
 {Modular Lattice for $C_{o}$-Operators.}

\author{Yun-Su Kim.  }

\address{Department of Mathematics, Indiana University, Bloomington,
Indiana, U.S.A. }

\email{kimys@indiana.edu}


\dedicatory{}



\newpage

\begin{abstract}We study modularity of the lattice \texttt{Lat}$(T)$ of
closed invariant subspaces for a $C_0$-operator $T$ and find a
condition such that \texttt{Lat}$(T)$ is a modular. Furthermore,
we provide a quasiaffinity preserving modularity.
\end{abstract}

\maketitle

\section*{Introduction}
A partially ordered set is said to be a \emph{lattice} if any two
elements \textbf{M} and \textbf{N} of it have a least upper bound
or supremum denoted by $\textbf{M}\vee{\textbf{N}}$ and a greatest
lower bound or infimum denoted by $\textbf{M}\cap\textbf{N}$. For
a Hilbert space $H$, $L(H)$ denotes the set of all bounded linear
operators from $H$ into $H$. For an operator $T$ in $L(H)$, the
set \texttt{Lat}$(T)$ of all closed invariant subspaces for $T$ is
a lattice. For \textbf{L}, \textbf{M}, and \textbf{N} in
\texttt{Lat}$(T)$ such that $\textbf{N}\subset{\textbf{L}}$, if
following identity is satisfied :
\begin{center}$\textbf{L}\cap(\textbf{M}\vee{\textbf{N}})=(\textbf{L}\cap{\textbf{M}})\vee{\textbf{N}},$
\end{center}
then \texttt{Lat}$(T)$ is called \emph{modular}. We study
\texttt{Lat}$(T)$ where $T$ is a $C_{0}$-operator which were first
studied in detail by B.Sz.-Nagy and C. Foias \cite{S2}. In this
paper \textbf{D} denotes the open unit disk in the complex
plane.\vskip0.1cm This paper is organized as follows. Section 1
contains preliminaries about operators of class $C_{0}$ and the
Jordan model of $C_0$-operators.

For operators $T_{1}\in{L(H_{1})}$ and $T_{2}\in{L(H_{2})}$, if
$X\in\{A\in{L(H):AT_{1}=T_{2}A}\}$, then we define a function
$X_{\ast}:\texttt{Lat}(T_{1})\rightarrow\texttt{Lat}(T_{2})$ as
following:
\begin{center}$X_{\ast}(M)=(XM)^{-}.$\end{center}

In Theorem \ref{x3}, we provide a quasiaffinity $Y$ such that
$Y_{\ast}$ preserves modularity. Furthermore, in section 2, we
provide a definition and prove some fundamental results of
\emph{property (P)} which was introduced by H. Bercovici
\cite{B4}.

In Theorem \ref{97}, we prove that if $T\in{L(H)}$ is an operator
of class $C_{0}$ with property $(P)$, then \texttt{Lat}$(T)$ is a
modular lattice.

The author would like to express her gratitude to her thesis
advisor, Professor Hari Bercovici.\vskip5cm

\newpage


\section{\textbf{$C_{0}$-Operators Relative to \textbf{D}}}\label{13}
\subsection{\textbf{A Functional Calculus.}}\label{11}
It is well-known that for every linear operator $A$ on a finite
dimensional vector space $V$ over the field $F$, there is a
minimal polynomial for $A$ which is the (unique) monic generator
of the ideal of polynomials over $F$ which annihilate $A$. If the
dimension of $F$ is not finite, then generally there is no such a
polynomial. However, to provide a function similar to a minimal
polynomial, B. Sz.-Nagy and C. Foias focused on a contraction
$T\in{L(H)}$ which is called to be \emph{completely nonunitary},
i.e. there is no invariant subspace $M$ for $T$ such that the
restriction $T|{M}$ of $T$ to the space $M$ is a unitary operator.

Let $H$ be a subspace of a Hilbert space $K$ and $P_{H}$ be the
orthogonal projection from $K$ onto $H$. We recall that if
$A\in{L(K)}$, and $T\in{L(H)}$, then $A$ is said to be a
\emph{dilation} of $T$ provided that for $n=1, 2, ...$,

\begin{equation}T^{n}=P_{H}A^{n}|H.\end{equation}

If $A$ is an isometry (unitary operator) then $A$ will be called
an \emph{isometric} (\emph{unitary}) \emph{dilation} of $T$. An
isometric (unitary) dilation $A$ of $T$ is said to be
\emph{minimal} if no restriction of $A$ to an invariant subspace
is an isometric (unitary) dilation of $T$. B. Sz.-Nagy proved the
following interesting result: \begin{prop}\cite{S2}\label{77}
Every contraction has a unitary dilation.\end{prop}

Let $T\in{L(H)}$ be a completely nonunitary contraction with
minimal unitary dilation $U\in{L(K)}$. For every polynomial
$p(z)=\sum_{j=0}^{n}a_{j}z^{j}$ we have
\begin{equation}\label{78}
p(T)=P_{H}p(U)|H,\end{equation} and so this formula suggests that
the functional calculus $p\rightarrow{p(T)}$ might be extended to
more general functions $p$. Since the mapping $p\rightarrow{p(T)}$
is a homomorphism from the algebra of polynomials to the algebra
of operators, we will extend it to a mapping which is also a
homomorphism from an algebra to the algebra of operators. By
Spectral Theorem, since $U\in{L(H)}$ is a normal operator, there
is a unique \emph{spectral measure} $E$ on the Borel subsets of
the spectrum of $U$ denoted as usual by $\sigma(U)$ such that
\begin{equation}\label{81} U=\int_{\sigma(U)}zdE(z).\end{equation}
Since the spectral measure $E$ of $U$ is absolutely continuous
with respect to Lebesgue measure on $\partial{\textbf{D}}$, for
$g\in{L^{\infty}(\sigma(U),E)}$, $g(U)$ can be defined as
follows:\begin{equation}\label{80}
g(U)=\int_{\sigma(U)}g(z)dE(z).\end{equation} It is clear that if
$g$ is a polynomial, then this definition agrees with the
preceding one. Since the spectral measure of $U$ is absolutely
continuous with respect to Lebesgue measure on
$\partial{\textbf{D}}$, the expression $g(U)$ makes sense for
every $g\in{L^{\infty}}=L^{\infty}(\partial{\textbf{D}})$. We
generalize formula (\ref{78}), and so for $g\in{L^{\infty}}$,
define $g(T)$ by
\begin{equation}\label{82}g(T)=P_{H}g(U)|H
    .\end{equation}

While the mapping $g\rightarrow{g(T)}$ is obviously linear, it is
not generally multiplicative, i.e. it is not a homomorphism.
Evidently it is convenient to find a subalgebra in $L^{\infty}$ on
which the functional calculus is multiplicative. Recall that
$H^{\infty}$ is the Banach space of all (complex-valued) bounded
analytic functions on the open unit disk $\textbf{D}$ with
supremum norm \cite{S2}. It turns out that $H^{\infty}$ is the
unique maximal algebra making the map a homomorphism between
algebras. We know that $H^{\infty}$ can be regarded as a
subalgebra of $L^{\infty}(\partial{\textbf{D}})$ \cite{B1}.

We note that the functional calculus with $H^{\infty}$ functions
can be defined in terms of independent of the minimal unitary
dilation. Indeed, if $u(z)=\sum_{n=0}^{\infty}a_{n}z^{n}$ is in
$H^{\infty}$, then
\begin{equation}\label{83}
u(T)=\lim_{r\rightarrow{1}}u(rT)=\lim_{r\rightarrow{1}}\sum_{n=0}^{\infty}a_{n}r^{n}T^{n},
\end{equation}
where the limit exists in the strong operator topology.

B. Sz.-Nagy and C. Foias introduced this important functional
calculus for completely nonunitary contractions.
\begin{prop}\label{12}Let $T\in{L(H)}$ be a completely nonunitary
contraction. Then there is a unique algebra representation
$\Phi_{T}$ from $H^{\infty}$ into $L(H)$ such that :\vskip0.2cm

(i) $\Phi_{T}(1)=I_{H}$, where $I_{H}\in{L(H)}$ is the identity
operator;

(ii) $\Phi_{T}(g)=T$, if $g(z)=z$ for all $z\in\textbf{D}$;

(iii) $\Phi_{T}$ is continuous when $H^{\infty}$ and $L(H)$ are
given the weak$^\ast$-topology.

(iv) $\Phi_{T}$ is contractive, i.e.
$\norm{\Phi_{T}(u)}\leq\norm{u}$ for all
$u\in{H^{\infty}}$.\end{prop}

We simply denote by $u(T)$ the operator $\Phi_{T}(u)$.

B.Sz.- Nagy and C. Foias \cite{S2} defined the \emph{class
$C_{0}$} relative to the open unit disk \textbf{D} consisting of
completely nonunitary contractions $T$ on $H$ such that the kernel
of $\Phi_{T}$ is not trivial.  If $T\in{L(H)}$ is an operator of
class $C_{0}$, then \begin{center}ker
$\Phi_{T}=\{u\in{H^{\infty}}:u(T)=0\}$\end{center} is a
weak$^{*}$-closed ideal of $H^{\infty}$, and hence there is an
inner function generating ker $\Phi_{T}$. The \emph{minimal
function} $m_{T}$ of an operator of class $C_{0}$ is the generator
of ker $\Phi_{T}$, and it seems as a substitute for the minimal
polynomial. Also, $m_{T}$ is uniquely determined up to a constant
scalar factor of absolute value one \cite{B1}. The theory of class
$C_{0}$ relative to the open unit disk has been developed by
B.Sz.- Nagy, C. Foias (\cite{S2}) and H. Bercovici (\cite{B1}).

\subsection{\textbf{Jordan Operator}} We know that every $n\times{n}$ matrix over an
algebraically closed field $F$ is similar to a unique Jordan
canonical form. To extend that theory to the $C_{0}$ operator
$T\in{L(H)}$, B.Sz.- Nagy and C. Foias \cite{S2} introduced a
weaker notion of equivalence. They defined a \emph{quasiaffine
transform} of $T$ which is bounded operator $T^\prime$ defined on
a Hilbert space $H^\prime$ such that there exists an injective
operator $X\in{L(H,H^{\prime})}$ with dense range in $H^{\prime}$
satisfying $T^{\prime}X=XT$. We write
\begin{center}$T\prec{T^{\prime}}$\end{center} if $T$ is a quasiaffine transform of
$T^{\prime}$. Instead of similarity, they introduced
\emph{quasisimilarity} of two operators, namely, $T$ and
$T^{\prime}$ are \emph{quasisimilar}, denoted by
\begin{center}$T\sim{T^{\prime}}$,\end{center} if $T\prec{T^{\prime}}$ and
$T^{\prime}\prec{T}$. \vskip0.2cm

Given an inner function $\theta\in{H^{\infty}}$, the \emph{Jordan
block }$S(\theta)$ is the operator acting on
$H(\theta)=H^{2}\ominus\theta{H^{2}}$, which means the orthogonal
complement of $\theta{H^2}$ in the Hardy space $H^2$, as follows
:\begin{equation} S(\theta)=P_{H(\theta)}S|H(\theta)\end{equation}
where $S\in{L(H^{2}})$ is the unilateral shift operator defined by
\begin{center}$(Sf)(z)=zf(z)$\end{center} and $P_{H(\theta)}\in{L(H^{2}})$ denotes the
orthogonal projection of $H^2$ onto ${H(\theta)}$.
\begin{prop}\cite{B1}\label{14}
For every inner function $\theta$ in $H^{\infty}$, the operator
$S(\theta)$ is of class $C_{0}$ and its minimal function is
$\theta$.\end{prop}

Let $\theta$ and $\theta^{\prime}$ be two inner functions in
$H^{\infty}$. We say that $\theta$ \emph{divides}
$\theta^{\prime}$ $($or $\theta$$\mid$$\theta^{\prime}$$)$ if
$\theta^{\prime}$ can be written as
$\theta^{\prime}=\theta$$\cdot$$\phi$ for some
$\phi\in{H^{\infty}}$. It is clear that $\phi\in{H^{\infty}}$ is
also inner. We will use the notation
$\theta\equiv{\theta}^{\prime}$ if
$\theta$$\mid$${\theta}^{\prime}$ and
${\theta}^{\prime}$$\mid$${\theta}$.

\begin{prop}\cite{B1}\label{z}
Let $T_1\in{L(H)}$ and $T_2\in{L(H)}$ be two completely nonunitary
contarctions of class $C_0$. If $T_1$ and $T_2$ are quasisimilar,
then $m_{T_1}\equiv{m_{T_2}}$.
\end{prop}
From Proposition \ref{14} and Proposition \ref{z}, we can easily
see that for every inner functions $\theta_1$ and $\theta_2$ in
$H^{\infty}$, if $S(\theta_{1})$ and $S(\theta_{2})$ are
quasisimilar, then $\theta_{1}\equiv\theta_{2}$. Conversely,
\begin{prop}\cite{B1}\label{03}
Let $\theta_1$ and $\theta_2$ be inner functions in $H^{\infty}$.
If $\theta_{1}\equiv\theta_{2}$, then $S(\theta_{1})$ and
$S(\theta_{2})$ are quasisimilar.
\end{prop}

 Let $\gamma$
be a cardinal number and
\begin{center}$\Theta=\{\theta_{\alpha}\in{H^{\infty}}:\alpha<\gamma\}$\end{center} be a
family of inner functions. Then $\Theta$ is called a \emph{model
function} if $\theta_{\alpha}|\theta_{\beta}$ whenever
card$(\beta)\leq$card$(\alpha)<\gamma$. The \emph{Jordan operator}
$S(\Theta)$ determined by the model function $\Theta$ is the
$C_{0}$ operator defined as
\begin{center}$S(\Theta)=\bigoplus_{\alpha<\gamma^{\prime}}S(\theta_{\alpha})$\end{center}
where $\gamma^{\prime}=$ min$\{\beta: \theta_{\beta}\equiv1\}$.

 We will call $S(\Theta)$ the \emph{Jordan model }of the operator
$T$ if \begin{center}$S(\Theta)\sim{T}$,\end{center} and in the
sequel $\bigoplus$$_{i<\gamma^{\prime}}$$S(\theta_{i}) $ always
means a \emph{Jordan operator} determined by a model
function.\vskip0.2cm By using Jordan blocks, $C_{0}$-operators
relative to the open unit disk \textbf{D} can be classified
(\cite{B1} Theorem 5.1) :\vskip0.2cm
\begin{thm}\label{15}
Any $C_{0}$-operator $T$ relative to the open unit disk \textbf{D}
acting on a Hilbert space is quasisimilar to a unique Jordan
operator.
\end{thm}
\vskip0.2cm
\begin{thm}\label{a4}
If $\Theta$ and $\Theta^{\prime}$ are two model functions and
$S(\Theta)\prec{S(\Theta^{\prime})}$, then
$\Theta\equiv\Theta^{\prime}$ and hence
$S(\Theta)=S(\Theta^{\prime})$.
\end{thm}
From Theorem \ref{15} and Theorem \ref{a4}, we can conclude that
$"\prec"$ is an equivalence relation on the set of
$C_{0}$-operators.

\section{\textbf{Lattice of subspaces}}\label{72}
\subsection{\textbf{Modular Lattice.}}Let $H$ be a Hilbert space. If $F_{i}(i\in{I})$ is a
subset of $H$, then the closed linear span of $\bigcup_{i}F_{i}$
will be denoted by $\bigvee_{i}F_{i}$. The collection of all
subspaces of a Hilbert space is a \emph{lattice}. This means that
the collection is partially ordered (by inclusion), and that any
two elements \textbf{M} and \textbf{N} of it have a least upper
bound or supremum (namely the span $\textbf{M}\vee{\textbf{N}}$)
and a greatest lower bound or infimum (namely the intersection
$\textbf{M}\cap{\textbf{N}}$). A lattice is called
\emph{distributive} if
\begin{equation}\label{71}\textbf{L}\cap(\textbf{M}\vee{\textbf{N}})=(\textbf{L}\cap{\textbf{M}})\vee(\textbf{L}\cap{\textbf{N}})
\end{equation}
for any element \textbf{L}, \textbf{M}, and \textbf{N} in the
lattice.

In the equation (\ref{71}), if $\textbf{N}\subset{\textbf{L}}$,
then $\textbf{L}\cap{\textbf{N}}=\textbf{N}$, and so the identity
becomes
\begin{equation}\label{72}\textbf{L}\cap(\textbf{M}\vee{\textbf{N}})=(\textbf{L}\cap{\textbf{M}})\vee{\textbf{N}}
\end{equation}

If the identity (\ref{72}) is satisfied whenever
$\textbf{N}\subset{\textbf{L}}$, then the lattice is called
\emph{modular}.

For an arbitrary operator $T\in{L(H)}$, \texttt{Lat}($T$) denotes
the collection of all closed invariant subspaces for $T$. The
following fact is well-known \cite{H0}.

\begin{prop}\label{74}
The lattice of subspaces of a Hilbert space H is modular if and
only if $\dim$ H is finite.
\end{prop}

We will think about \texttt{Lat}$(T)$ for a $C_{0}$-operator $T$.

\begin{defn}\label{75}
The \emph{cyclic multiplicity} $\mu_{T}$ of an operator
$T\in{L(H)}$ is the smallest cardinal of a subset $A\subset{H}$
with the property that $\bigvee_{n=0}^{\infty}T^{n}A=H$. The
operator $T$ is said to be \emph{multiplicity-free} if
$\mu_{T}=1$.
\end{defn}
Thus $\mu_{T}$ is the smallest number of cyclic subspaces for $T$
that are needed to generate $H$, and $T$ is multiplicity-free if
and only if it has a cyclic vector.

\subsection{\textbf{Property $(P)$.}}\label{73}

Let $H$ be a Hilbert space and for an operator $T\in{L(H)}$,
$T^{\ast}$ denote the adjoint of $T$. It is well known that $H$ is
finite-dimensional if and only if every operator $X\in{L(H)}$,
with the property $\ker (X)=\{0\}$, also satisfies
$\ker(X^{\ast})=\{0\}$. The following definition is a natural
extension of finite dimensionality.
\begin{defn}\label{94}
An operator $T\in{L(H)}$ is said to have \emph{property (P) } if
every operator $X\in\{T\}^{\prime}$ with the property that
$\ker(X)=\{0\}$ is a quasiaffinity, i.e.,
$\ker(X^{\ast})=\ker(X)=\{0\}$.
\end{defn}
From the fact that the commutant $\{0\}^{\prime}$ of zero operator
on $H$ coincides with $L(H)$, we can see that $H$ is
finite-dimensional if and only if the zero operator on $H$ has
property $(P)$.

Let $T_{1}$ and $T_{2}$ be operators in $L(H)$. Suppose that
\begin{center}$X\in\{A\in{L(H):AT_{1}=T_{2}A}\}$.\end{center} If $M$ is in
\texttt{Lat}$(T_{1})$, then $(XM)^-$ is in \texttt{Lat}$(T_{2})$.
By using these facts, we define a function $X_{\ast}$ from
\texttt{Lat}$(T_{1})$ to \texttt{Lat}$(T_{2})$ as following :
\begin{equation}\label{95}X_{\ast}(M)=(XM)^{-}.\end{equation}
The operator $X$ is said to be a
$(T_{1},T_{2})$\emph{-lattice-isomorphism} if $X_{\ast}$ is a
bijection of \texttt{Lat}$(T_{1})$ onto \texttt{Lat}$(T_{2})$. We
will use the name lattice-isomorphism instead of
$(T_{1},T_{2})$-lattice-isomorphism if no confusion may arise.

If $X\in\{A\in{L(H):AT_{1}=T_{2}A}\}$, then
$X^{\ast}T_{2}^{\ast}=T_{1}^{\ast}X^{\ast}$. Thus
$(X^{\ast})_{\ast}:\texttt{Lat}(T_{2}^{\ast})\rightarrow\texttt{Lat}(T_{1}^{\ast})$
is well-defined by
\begin{center}$(X^{\ast})_{\ast}(M^{\prime})=(X^{\ast}M^{\prime})^{-}$.\end{center}
\begin{prop}\label{a2}\cite{B1} (\texttt{Theorem 7.1.9})
Suppose that $T\in{L(H)}$ is an operator of class $C_{0}$ with
Jordan model $\bigoplus_{\alpha}S(\theta_{\alpha})$. Then $T$ has
property $(P)$ if and only if
\begin{center}$\bigwedge_{j<\omega}\theta_{j}\equiv{1}.$\end{center}
Thus, if $T$ has property $(P)$, then $H$ is separable and
$T^{\ast}$ also has property $(P)$.
\end{prop}
\begin{prop}\label{t}\cite{B1}
An operator $T$ of class $C_{0}$ fails to have property $(P)$ if
and only if $T$ is quasisimilar to $T|N$, where $N$ is a proper
invariant subspace for $T$.
\end{prop}

\begin{prop}\label{02}\cite{B1}(\texttt{Lemma 7.1.20})
Assume that $T_{1}\in{L(H_1)}$ and $T_{2}\in{L(H_2)}$ are two
operators, and $X\in\{A\in{L(H_{1},H_{2}):AT_{1}=T_{2}A}\}$. If
the mapping $X_{\ast}$ is onto \texttt{Lat}$(T_{2})$ if and only
if $(X^{\ast})_{\ast}$ is one-to-one on
\texttt{Lat}$({T_{2}^{\ast}})$.
\end{prop}

\begin{cor}\label{s}
Assume that $T_{1}\in{L(H_1)}$ and $T_{2}\in{L(H_2)}$ are two
operators, and $X\in\{A\in{L(H_{1},H_{2}):AT_{1}=T_{2}A}\}$. The
mapping $X_{\ast}$ is one-to-one on \texttt{Lat}$(T_{1})$ if and
only if $(X^{\ast})_{\ast}$ is onto
\texttt{Lat}$({T_{1}^{\ast}})$.
\end{cor}
\begin{proof}
Since $XT_{1}=T_{2}X$,
$T_{1}^{\ast}X^{\ast}=X^{\ast}T_{2}^{\ast}$. By Proposition
\ref{02}, $(X^{\ast})_{\ast}$ is onto \texttt{Lat}$(T_{1}^{\ast})$
if and only if $(X^{\ast\ast})_{\ast}=X_{\ast}$ is one-to-one on
\texttt{Lat}$(T_{1})$.
\end{proof}

From Proposition \ref{02} and Corollary \ref{s}, we obtain the
following result.
\begin{cor}\label{s3}
If $T_{1}\in{L(H_1)}$ and $T_{2}\in{L(H_2)}$ are two operators,
and $X\in\{A\in{L(H_{1},H_{2}):AT_{1}=T_{2}A}\}$, then $X$ is a
lattice-isomorphism if and only if $X^{\ast}$ is a
lattice-isomorphism.
\end{cor}

\begin{prop}\label{a1}\cite{B1} (\texttt{Proposition 7.1.21})
Assume that $T_{1}\in{L(H_{1})}$ and $T_{2}\in{L(H_{2})}$ are two
quasisimilar operators of class $C_{0}$, and
$X\in\{A\in{L(H_{1},H_{2}):AT_{1}}$ $=T_{2}A\}$ is an injection.
If $T_{1}$ has property (P), then $X$ is a lattice-isomorphism.
\end{prop}

Recall that if $T$ is an operator on a Hilbert space, then ker
$T=$ $($ran ${T}^{\ast})^{\perp}$ and ker $T^\ast=$ $($ran
$T)^\perp$.

\begin{cor}\label{s1}
Assume that $T_{1}\in{L(H_{1})}$ and $T_{2}\in{L(H_{2})}$ are two
quasisimilar operators of class $C_{0}$, and
$X\in\{A\in{L(H_{1},H_{2}):AT_{1}=T_{2}A}\}$ has dense range. If
$T_{2}$ has property (P), then $X$ is a lattice-isomorphism.
\end{cor}
\begin{proof}
Since $XT_{1}=T_{2}X$,
$T_{1}^{\ast}X^{\ast}=X^{\ast}T_{2}^{\ast}$. Let $Y=X^{\ast}$ and
so \begin{equation}\label{s2} YT_{2}^{\ast}=T_{1}^{\ast}Y.
\end{equation}
From the fact that $\ker{Y}=\ker(X^{\ast})=( \texttt{ran }
{X})^{\perp}=\{0\}$, we conclude that $Y$ is injective. Since
$T_{2}$ has property $(P)$, so does $T_{2}^{\ast}$ by Proposition
\ref{a2}.
 By Proposition \ref{a1} and
equation (\ref{s2}), $Y=X^{\ast}$ ia a lattice-isomorphism. From
Corollary \ref{s3}, it is proven that $X$ is a
lattice-isomorphism.
\end{proof}

\begin{cor}\label{s4}
Suppose that $T_{i}\in{L(H_{i})}(i=1,2)$ is a $C_{0}$-operator and
$T_{1}$ has property $(P)$. If
$X\in\{A\in{L(H_{1},H_{2}):AT_{1}=T_{2}A}\}$ and $X$ is an
injection, then $X$ is a lattice-isomorphism.
\end{cor}
\begin{proof}
Define $Y:H_{1}\rightarrow(XH_{1})^{-}$ by \begin{center}$Yh=Xh$
for any $h\in{H_{1}}$.\end{center} Since $X$ is an injection, so
is $Y$. Clearly, $Y$ has dense range. Note that $(XH_{1})^{-}$ is
invariant for $T_{2}$. By definition of $Y$,
\begin{equation}\label{64}
YT_{1}=(T_{2}|(XH_{1})^{-})Y.
\end{equation}
It follows that $T_{1}\prec(T_{2}|(XH_{1})^{-})$ and so
$T_{1}\sim(T_{2}|(XH_{1})^{-})$. By Proposition \ref{a1}, it is
proven.
\end{proof}
\begin{cor}\label{t1}
Suppose that $T_{i}\in{L(H_{i})}(i=1,2)$ is a $C_{0}$-operator and
$T_{2}$ has property $(P)$. If
$X\in\{A\in{L(H_{1},H_{2}):AT_{1}=T_{2}A}\}$ and $X$ has a dense
range, then $X$ is a lattice-isomorphism.
\end{cor}
\begin{proof}
By assumption, $X^{\ast}T_{2}^{\ast}=T_{1}^{\ast}X^{\ast}$. Since
$T_{2}$ has property $(P)$, by Proposition \ref{a2}, so does
$T_{2}^{\ast}$.

Because $X$ has dense range, $X^{\ast}:H_{2}\rightarrow{H_{1}}$ is
an injection. By Corollary \ref{s4}, $X^{\ast}$ is a lattice
isomorphism. From Corollary \ref{s3}, $X$ is also a lattice
isomorphism.
\end{proof}

\subsection{\textbf{Quasi-Affinity and Modular Lattice}}\label{q6}
For operators $T_{1}\in{L(H_{1})}$ and $T_{2}\in{L(H_{2})}$, if
$Y\in\{B\in{L(H_{1},H_{2})}:BT_{1}=T_{2}B\}$, then we define a
function
\begin{center}$Y_{\ast}:{\texttt{Lat}(T_{1})}\rightarrow{\texttt{Lat}(T_{2})}$\end{center}
the same way as equation (\ref{95}). For any
$N\in{\texttt{Lat}(T_{2})}$, if $M=Y^{-1}(N)$, then
$YT_{1}(M)=T_{2}Y(M)\subset{T_{2}N}\subset{N}$ and so
$T_{1}(M)\subset{M}$. It follows that
\begin{center}$M=Y^{-1}(N)\in{\texttt{Lat}(T_{1})}$\end{center} for any
$N\in{\texttt{Lat}(T_{2})}$. If $Y$ is invertible, that is, $T_1$
and $T_{2}$ are similar, and $\texttt{Lat}(T_{1})$ is modular,
then clearly, $\texttt{Lat}(T_{2})$ is also modular. In this
section, we consider when $T_1$ and $T_2$ are quasi-similar
instead of similar, and find an assumption in Theorem \ref{x3}
such that $\texttt{Lat}(T_{2})$ is modular, whenever
$\texttt{Lat}(T_{1})$ is modular.
\begin{prop}\label{o}
Let $T_{1}\in{L(H_{1})}$ and $T_{2}\in{L(H_{2})}$. Suppose that
$Y\in\{B\in{L(H_{1},H_{2})}:BT_{1}=T_{2}B\}$ and for any
$N\in{\texttt{Lat}(T_{2})}$, the condition $M=Y^{-1}(N)$ implies
that $Y_{\ast}(M)=N$.

Then for any $M_{i}=Y^{-1}(N_i)$ with
$N_{i}\in{\texttt{Lat}(T_{2})}$ $(i=1,2)$,
\begin{center}$Y_{\ast}(M_{1}\cap{M_{2}})=Y_{\ast}(M_1)\cap{Y_{\ast}(M_2)}$.\end{center}
\end{prop}
\begin{proof}
Assume that $N_{i}\in{\texttt{Lat}(T_{2})}$ and
$M_{i}=Y^{-1}(N_i)$ for $i=1,2$. Then by assumption, we obtain
\begin{equation}\label{o1}Y_{\ast}(M_{i})=N_{i}.\end{equation} Since
$Y^{-1}(N_{1}\cap{N_{2}})=Y^{-1}(N_{1})\cap{Y}^{-1}(N_2)=M_{1}\cap{M_2},$
by assumption,
\begin{center}$Y_{\ast}(M_{1}\cap{M_2})=N_{1}\cap{N_2}$\end{center} which
proves that
$Y_{\ast}(M_{1}\cap{M_2})=Y_{\ast}(M_1)\cap{Y_{\ast}(M_2)}$ by
equation (\ref{o1}).

\end{proof}

\begin{thm}\label{x3}
Let $T_{1}\in{L(H_{1})}$ be a quasiaffine transform of
$T_{2}\in{L(H_{2})}$ and
$Y\in\{B\in{L(H_{1},H_{2})}:BT_{1}=T_{2}B\}$ be a quasiaffinity.

If $Y_{\ast}:\texttt{Lat}(T_{1})\rightarrow\texttt{Lat}(T_{2})$ is
onto and $\texttt{Lat}(T_{1})$ is modular, then
$\texttt{Lat}(T_{2})$ is also modular.
\end{thm}
\begin{proof}
Suppose that $\texttt{Lat}(T_{2})$ is not modular. Then there are
invariant subspaces $N_{i}(i=1,2,3)$ for $T_{2}$ such that
\begin{equation}\label{f}N_{3}\subset{N_{1}},\end{equation} and
\begin{center}$(N_{1}\cap{N_2})\vee{N_3}\neq{N}_{1}\cap(N_{2}\vee{N_3}).$\end{center}

Let \begin{equation}\label{f2}M_{i}=Y^{-1}(N_i),\end{equation} for
$i=1,2,3.$ Since $YT_{1}=T_{2}Y$, definition (\ref{f2}) of $M_{i}$
implies that for $i=1,2,3,$
\begin{center}$YT_{1}(M_{i})=T_{2}Y(M_{i})\subset{T_{2}N_{i}}\subset{N_{i}}.$\end{center}
It follows that $T_{1}M_{i}\subset{Y^{-1}(N_i)}=M_{i}$ for
$i=1,2,3$. Thus $M_{i}$ is a closed invariant subspace for
$T_{1}$. Condition (\ref{f}) implies that
\begin{center}\label{f4}$M_{3}\subset{M_1}.$\end{center}

Since $Y(M_i)\subset{N_i}$, for $i=1,2,3,$
\begin{equation}\label{f5}Y_{\ast}(M_i)=(Y(M_i))^{-}\subset{N_i}.\end{equation}

Since $Y_{\ast}$ is onto, there is a function
$\phi:\texttt{Lat}(T_{2})\rightarrow\texttt{Lat}(T_{1})$ such that
$Y_{\ast}\circ\phi$ is the identity mapping on
\texttt{Lat}$(T_{2})$. Hence for $i=1,2,3,$
\begin{center}\label{f7}$
Y_{\ast}(\phi({N_i}))=Y(\phi({N_i}))^{-}=N_{i}.$\end{center} It
follows that for $i=1,2,3,$\begin{equation}\label{f8}
\phi({N_i})\subset{M_{i}}.\end{equation} Since $Y_{\ast}\circ\phi$
is the identity mapping on \texttt{Lat}$(T_{2})$, (\ref{f8})
implies that for $i=1,2,3,$
\begin{equation}\label{f9}N_{i}=Y_{\ast}(\phi({N_i}))\subset{Y}_{\ast}(M_{i}).\end{equation}
By (\ref{f5}) and (\ref{f9}), we get \begin{equation}\label{g}
Y_{\ast}(M_{i})=N_{i},\end{equation} for $i=1,2,3.$ Hence we can
easily see that function $Y$ satisfies the assumptions of
Proposition \ref{o}.

Thus by Proposition \ref{o} and equation (\ref{g}),
\begin{equation}\label{o2}
Y_{\ast}[M_{1}\cap({M_2}\vee{M_3})]=Y_{\ast}(M_1)\cap{Y_{\ast}(M_{2}\vee{M_3})}=N_{1}\cap({N_2}\vee{N_3}).
\end{equation}

Since $M_{1}\cap{M_{2}}=Y^{-1}(N_{1})\cap{Y^{-1}}(N_{2})
=Y^{-1}(N_{1}\cap{N_{2}})$, by the same way as above, we obtain
\begin{equation}\label{k5}Y_{\ast}(M_{1}\cap{M_{2}})=N_{1}\cap{N_{2}}.\end{equation}
By equations (\ref{g}) and (\ref{k5}), we obtain
\begin{equation}\label{g1}Y_{\ast}[(M_{1}\cap{M_2})\vee{M_3}]=(N_{1}\cap{N_2})\vee{N_3}.\end{equation}

Since $(N_{1}\cap{N_2})\vee{N_3}\neq{N}_{1}\cap(N_{2}\vee{N_3})$,
from equations (\ref{o2}) and (\ref{g1}), we can conclude that
\begin{center}\label{g3}$(M_{1}\cap{M_2})\vee{M_3}\neq{M}_{1}\cap(M_{2}\vee{M_3}).$\end{center}
Therefore $\texttt{Lat}(T_{1})$ is  not modular.

\end{proof}

\section{\textbf{Modular Lattice for
$C_{0}$-Operators with Property $(P)$}}\label{q1} We provide some
operators, say $T$, of class $C_0$ such that $\texttt{Lat}(T)$ is
modular.
\begin{prop}\label{q}\cite{B1}
Let $\theta$ be a nonconstant inner function in $H^\infty$. Then
every invariant subspace $M$ of $S(\theta)$ has the form
\begin{center}$\phi{H^{2}}\ominus\theta{H^{2}}$\end{center} for some inner devisor
$\phi$ of $\theta$.
\end{prop}
We can easily check that if
$\textbf{M}_{1}=\theta_{1}{H^{2}}\ominus\theta{H^{2}}$ and
$\textbf{M}_{2}=\theta_{2}{H^{2}}\ominus\theta{H^{2}}$ where
$\theta_i$ ($i=1,2$) is an inner inner devisor of $\theta$, then
\begin{equation}\label{q2}\textbf{M}_{1}\cap\textbf{M}_{2}=(\theta_{1}\vee\theta_{2}){H^{2}}\ominus\theta{H^{2}}\end{equation}
and
\begin{equation}\label{q3}\textbf{M}_{1}\vee\textbf{M}_{2}=(\theta_{1}\wedge\theta_{2}){H^{2}}\ominus\theta{H^{2}}\end{equation}
where $\theta_{1}\wedge\theta_{2}$ and $\theta_{1}\vee\theta_{2}$
denote the greatest common inner divisor and least common inner
multiple of $\theta_1$ and $\theta_2$, respectively. Note that if
$\textbf{M}_{1}\subset\textbf{M}_{2}$, then
\begin{equation}\label{q7}\theta_{2}|\theta_{1}.\end{equation}

\begin{lem}
If $\theta$ is an inner function in $H^\infty$, then
$\texttt{Lat}(S(\theta))$ is distributive.
\end{lem}
\begin{proof}
Let \textbf{M$_{1}$}, \textbf{M$_2$}, and \textbf{M$_3$} be
invariant subspaces for $S(\theta)$. Then by Proposition \ref{q},
there are nonconstant inner functions $\theta_1$, $\theta_2$, and
$\theta_3$ in $H^\infty$ such that
\begin{center}$\textbf{M}_{i}=\theta_{i}{H^{2}}\ominus\theta{H^{2}}$ for
$i=1,2,3.$\end{center} From equations (\ref{q2}) and (\ref{q3}),
we obtain that
\begin{equation}\label{q4}\textbf{M}_{1}\cap(\textbf{M}_{2}\vee{\textbf{M}_{3}})=
(\theta_{1}\vee(\theta_{2}\wedge\theta_{3})){H^{2}}\ominus\theta{H^{2}},\end{equation}
and
\begin{equation}\label{q5}(\textbf{M}_{1}\cap\textbf{M}_{2})\vee({\textbf{M}_{1}\cap\textbf{M}_{3}})=
((\theta_{1}\vee\theta_{2})\wedge(\theta_{1}\vee\theta_{3}){H^{2}}\ominus\theta{H^{2}}.\end{equation}
Since
$\theta_{1}\vee(\theta_{2}\wedge\theta_{3})=(\theta_{1}\vee\theta_{2})\wedge(\theta_{1}\vee\theta_{3})$,
by equations (\ref{q4}) and (\ref{q5}), this lemma is proven.

\end{proof}

 In this section, we
will consider a sufficient condition for \texttt{Lat}($T$) of a
$C_{0}$-operator $T$ to be modular.

\begin{prop}\label{t2}\cite{B1} (\texttt{Proposition 2.4.3})
Let $T\in{L(H)}$ be a completely nonunitary contraction, and $M$
be an invariant subspace for $T$. If
\begin{equation}T=\begin{pmatrix}
T_{1}&X\\0&T_{2}\end{pmatrix}\end{equation} is the
triangularization of $T$ with respect to the decomposition
$H=M\oplus(H$ $\ominus{M})$, then $T$ is of class $C_{0}$ if and
only if $T_1$ and $T_2$ are operators of class $C_{0}$.
\end{prop}
\begin{prop}\label{a3}\cite{B1} (\texttt{Corollary 7.1.17})
Let $T\in{L(H)}$ is an operator of class $C_{0}$, $M$ be an
invariant subspace for $T$, and \begin{equation}T=\begin{pmatrix}
T_{1}&X\\0&T_{2}\end{pmatrix}\end{equation} be the
triangularization of $T$ with respect to the decomposition
$H=M\oplus(H$ $\ominus{M})$. Then $T$ has property $(P)$ if and
only if $T_1$ and $T_2$ have property $(P)$.
\end{prop}

Let $H$ and $K$ be Hilbert spaces and $H\oplus{K}$ denote the
algebraic direct sum. Recall that $H\oplus{K}$ is also a Hilbert
space with an inner product

\centerline{$(\langle{h}_{1},k_{1}\rangle,\langle{h_2},k_{2}\rangle)
=(h_{1},h_{2})+(k_{1},k_{2})$}

\begin{thm}\label{97}
Let $T\in{L(H)}$ be an operator of class $C_{0}$ with property
$(P)$. Then \texttt{Lat}$(T)$ is a modular lattice.
\end{thm}
\begin{proof}
Suppose that $T$ has property $(P)$ and let $M_{i}$ $(i=1,2,3)$ be
an invariant subspace for $T$ such that $M_{3}\subset{M_{1}}$.
Then evidently,
\begin{equation}\label{07}(M_{1}\cap{M_{2}})\vee{M_3}\subset{M_{1}\cap(M_{2}\vee{M_3})}.\end{equation}

Let $T_{i}=T|M_{i}$ $(i=1,2,3)$. Define a linear transformation
$X:M_{2}\oplus{M_3}\rightarrow{M_{2}\vee{M_{3}}}$ by
\begin{center}$X(a_{2}\oplus{a_{3}})=a_{2}+a_{3}$\end{center} for $a_{2}\in{M_{2}}$ and
$a_{3}\in{M_3}$.\vskip0.1cm Then for
$a_{2}\oplus{a_{3}}\in{M_{2}\oplus{M_3}}$ with
$\norm{a_{2}\oplus{a_{3}}}\leq{1}$,
$\norm{X(a_{2}\oplus{a_{3}})}=\norm{a_{2}+a_{3}}\leq\norm{a_{2}}+\norm{a_{3}}\leq{2}$.
It follows that $\norm{X}\leq{2}$ and so $X$ is bounded.

 Since
$M_{2}\vee{M_{3}}$ is generated by $\{a_{2}+a_{3}:a_{2}\in{M_{2}}$
and $a_{3}\in{M_3}\}$, $X$ has dense range. By definition of
$T_{i}$ $(i=1,2,3)$,
\begin{center}$X(T_{2}\oplus{T_3})(a_{2}\oplus{a_{3}})=Ta_{2}+Ta_{3}$\end{center} and
\begin{center}$(T|M_{2}\vee{M_3})X(a_{2}\oplus{a_{3}})=Ta_{2}+Ta_{3}$.\end{center} Thus
\begin{center}$X(T_{2}\oplus{T_3})=(T|M_{2}\vee{M_3})X.$\end{center}

By Proposition \ref{t2}, $T_{2}\oplus{T_{3}}$ and
$T|M_{2}\vee{M_3}$ are of class $C_{0}$ and since $T$ has property
$(P)$, by Proposition \ref{a3}, we conclude that
$T|M_{2}\vee{M_3}$ also has Property $(P)$. By Corollary \ref{t1},
$X$ is a lattice-isomorphism.

Thus
$X_{\ast}:\texttt{Lat}(T_{2}\oplus{T_3})\rightarrow\texttt{Lat}(T|M_{2}\vee{M_3})
$ is onto. Let
\begin{equation}\label{99}M=\{a_{2}\oplus{a_{3}}\in{M_{2}\oplus{M_3}}:a_{2}+a_{3}\in{M_1}\}.\end{equation}
Since $M=X^{-1}(M_1)$, $M$ is a closed subspace of
$M_{2}\oplus{M_3}$. Evidently, $M$ is invariant for
$T_{2}\oplus{T_3}$. From the equation (\ref{99}), we conclude that
\begin{equation}\label{100}M=(M_{1}\cap{M_2})\oplus{M_3}.\end{equation}
Since $X^{-1}(M_{1}\cap(M_{2}\vee{M_{3}}))
=\{a_{2}\oplus{a_{3}}\in{M_{2}\oplus{M_3}}:a_{2}+a_{3}\in{M_{1}\cap(M_{2}\vee{M_{3}})}\}
=\{a_{2}\oplus{a_{3}}\in{M_{2}\oplus{M_3}}:a_{2}+a_{3}\in{M_{1}}\}$,

\begin{center}$X^{-1}(M_{1}\cap(M_{2}\vee{M_{3}}))=M$\end{center}
Since $X$ is a lattice-isomorphism,
\begin{equation}\label{04}X_{\ast}M=(XM)^{-}=M_{1}\cap(M_{2}\vee{M_{3}}).\end{equation}

By equation (\ref{100}) and definition of $X$,
\begin{equation}\label{06}X_{\ast}M=(XM)^{-}\subset{(M_{1}\cap{M_{2}})\vee{M_{3}}}.\end{equation}
From (\ref{04}) and (\ref{06}), we conclude that
\begin{equation}\label{08}
M_{1}\cap(M_{2}\vee{M_{3}})\subset{(M_{1}\cap{M_{2}})\vee{M_{3}}}.\end{equation}
Thus if $T$ has property $(P)$, then by (\ref{07}) and (\ref{08}),
we obtain that
\begin{center}$M_{1}\cap(M_{2}\vee{M_{3}})={(M_{1}\cap{M_{2}})\vee{M_{3}}}$.\end{center}
\vskip0.2cm

\end{proof}

------------------------------------------------------------------------

\bibliographystyle{amsplain}
\bibliography{xbib}
\end{document}